\documentclass[11pt]{smfart}

\usepackage{amsmath,a4wide}
\usepackage{amssymb}
\usepackage{xspace}
\usepackage{euscript}
\usepackage{mathrsfs} 
\usepackage[all]{xy}
\usepackage{hyperref}
\usepackage[francais]{babel}
\usepackage[T1]{fontenc}
\usepackage{textcomp}
\usepackage{smfthm}

\newtheorem{proposition}[section]{Proposition}	

\newtheorem{theoreme}[section]{Th\'eor\`eme}

\theoremstyle{definition}

\newenvironment{demo}{\noindent{\textit{D\'emonstration. --- }}}{~\qedsymbol \vspace{4mm}}

\renewcommand{\theenumi}{\textbf{\roman{enumi}}} 

\renewcommand\theequation{\thesection.\textbf{\alph{equation}}} 

\makeatletter
\@addtoreset{equation}{subsection}
\makeatother

\renewcommand\thesubsubsection{\textbf{\thesubsection.\arabic{subsubsection}}}
\renewcommand\thesubsection{\textbf{\thesection.\arabic{subsection}}}

\newcommand{\Q}{\mathbb Q}

\newcommand{\C}{\mathbb C}

\newcommand{\Spec}{\mathrm{Spec}}

\newcommand{\iso}{\buildrel{\sim}\over\rightarrow}

\renewcommand{\L}{\mathrm{L}}		

\def\commutatif{\ar@{}[rd]|{\circlearrowleft}}
\def\cartesien{\ar@{}[rd]|{\square}}

\renewcommand{\lim}{{\mathrm{lim}}} 

\def\egazéro#1#2{[{\bf \'EGA}~$0_{\textsc{#1}}$~#2]}		

\numberwithin{equation}{subsubsection}
\renewcommand{\theequation}{\thesubsection.\Alph{equation}}
\newenvironment{demoth1}{\noindent{\textit{D\'emonstration du théorème 1. --- }}}{~\qedsymbol \vspace{4mm}}

\title{Cohomologie relative des formes cuspidales}
\author{Beno\^it Stroh}
\date{30 juillet 2013}
\email{stroh@math.univ-paris13.fr}
\address{C.N.R.S, Université Paris 13, LAGA,
99 avenue J.B. Clément,
93430 Villetaneuse
France}

\begin{document}

\maketitle

Dans cette note, nous proposons une démonstration élémentaire du théorème suivant. Elle résulte d'une discussion vietnamienne avec Hélène Esnault, que nous remercions vivement.

Soit~$X_E$ une variété de Shimura définie sur son corps réflexe~$E$ qui est une extension finie de~$\Q$. D'après~\cite{Pi}, il existe une compactification minimale~$X^*_E$ et une compactification toroïdale~$\bar{X}_E$ définies sur~$E$. Le schéma~$X_E^*$ est par construction projectif muni d'un faisceau ample~$\omega_E$. Il existe par ailleurs un morphisme propre~$\pi_E : \bar{X}_E \rightarrow X^*_E$ qui est un isomorphisme sur~$X_E$. Quitte à bien choisir la combinatoire dont dépend la compactification toroïdale, on peut supposer que~$\bar{X}_E$ est lisse sur~$E$, que le bord réduit~$D_E$ de~$X_E$ dans~$\bar{X}_E$ est un diviseur à croisements normaux et que~$-D_E$ est relativement ample pour~$\pi_E$. On dispose enfin de l'isomorphisme de Kodaira-Spencer qui garantit l'existence d'un entier~$n$ tel que $\mathcal{K}_{\bar{X}_E}(\mathrm{log}(D_E))\iso \pi_E^*(\omega^{\otimes n})$ où~$\mathcal{K}_{\bar{X}_E}(\mathrm{log}(D))$ est le faisceau dualisant logarithmique de~$\bar{X}_E$.

Notons~$\mathcal{O}_E$ l'anneau des entiers de~$E$ et~$\Delta$ le discriminant de~$E$ sur~$\Q$. Soit~$\Spec(k)$ un point de~$\Spec(\mathcal{O}_E[1/\Delta])$. Lorsque~$k$ est de caractéristique nulle, on notera~$\mathcal{O}_k=k$ et lorsque~$k$ est de caractéristique positive, on notera~$\mathcal{O}_k=W(k)$ les vecteurs de Witt. On suppose que~$X_E$, $X^*_E$, $\bar{X}_E$, $D_E$ et~$\pi_E$ admettent un modèle $X$, $X^*$, $\bar{X}$, $D$ et~$\pi$ sur~$\Spec(\mathcal{O}_k)$. On demande que $X^*$ soit projectif et que~$\pi$ soit un isomorphisme sur~$X$. On demande aussi que~$\bar{X}$ soit lisse sur~$\Spec(\mathcal{O}_k)$, que le complémentaire réduit~$D$ de~$X$ dans~$\bar{X}$ soit un diviseur à croisements normaux relatif sur~$\Spec(\mathcal{O}_k)$ et que~$-D$ soit relativement ample pour~$\pi$. On demande enfin que l'isomorphisme de Kodaira-Spencer s'étende en un isomorphisme $\mathcal{K}_{\bar{X}}(\mathrm{log}(D))\iso \pi^*(\omega^{\otimes n})$.

Lorsque~$k$ est de caractéristique nulle, l'existence de tels modèles est tautologique. En caractéristique positive, de tels modèles ont été construits  dans~\cite{Lan} lorsque~$X$ est PEL et~$v$ est une place de bonne réduction pour~$X$.

\begin{theoreme} \label{thprinc} Si la caractéristique de~$k$ est non nulle, on suppose qu'elle est supérieure ou égale à la dimension de~$X_E$. On~a $\mathrm{R}^q \: \pi_* \mathcal{O}_{\bar{X}}(-D)=0$ sur~$X^*\times\Spec(k)$ pour tout~$q>0$.
\end{theoreme}

Ce théorème a été prouvé dans~\cite{AIP} pour les variétés de Siegel grâce à une analyse explicite des compactifications. Leur démonstration est plus compliquée que la n\^otre mais a l'avantage d'être valable en toute caractéristique. Ce théorème sera également démontré dans~\cite{HLTT} pour des variétés de Shimura PEL générales.
Il est utilisé dans plusieurs travaux récents (\cite{AIP}, \cite{ERX} et~\cite{PS} en caractéristique positive, \cite{HLTT} et~\cite{TX} en caractéristique nulle).

\medskip 

Commençons par prouver l'énoncé suivant, qui est une version faible du théorème d'annulation relative de Grauert-Riemenschneider valable en caractéristique non nulle~$p$.

\begin{proposition}\label{prop_generale} Soit $K$ un corps parfait de caractéristique~$p$ et $f:Z\rightarrow Y$ un morphisme birationnel entre schémas propres sur~$\Spec(W_2(K))$ où~$W_2(K)$ désigne les seconds vecteurs de Witt de~$K$. On suppose~$Y$ projectif et~$Z$ lisse sur~$\Spec(W_2(K))$. On suppose qu'il existe un diviseur à croisements normaux relatif~$D$ dans~$Z$ dont l'opposé est relativement ample pour~$f$. On suppose~$p\geq \mathrm{dim}(Z)$. Alors si~$\mathcal{K}_{Z}$ désigne le faisceau dualisant de~$X$, on a~$\mathrm{R}^q {f}_{1,*}\:  (\mathcal{K}_{Z}) = 0$ pour tout~$q>0$ où~${f}_1=f\times \Spec(K)$ est la réduction modulo~$p$ de~$f$. 
\end{proposition}

\begin{demo} On designe par~${Z}_1$, ${Y}_1$ et~${f}_1$ les réductions sur~$\Spec(K)$. Soit~$\mathcal{L}$ un faisceau ample sur~$Y$. Par hypothèse, $\mathcal{O}_X(-D)$ est relativement ample pour~$f$. D'après EGA II, prop.4.6.13.ii, il existe un entier~$n_0$ tel que pour tout~$n>n_0$ le faisceau~$f^* \mathcal{L}^n(-D)$ soit ample sur~$Z$. D'après la forme logarithmique du théorème d'annulation de Kodaira-Raynaud~\cite[2.8 et 4.2]{DI}, qui utilise l'existence d'un relèvement à $W_2(K)$ et l'hypothèse $p\geq \mathrm{dim}(Z)$, on a
$$\mathrm{H}^q({Z}_1,\: \mathcal{K}_Z(\mathrm{log}(D))\otimes f^* \mathcal{L}^n(-D)) \: = \: 0$$
pour tout~$q>0$. Donc
$$\mathrm{H}^q({Z}_1,\: \mathcal{K}_Z\otimes f^* \mathcal{L}^n) \: = \: 0$$
pour tout~$q>0$. Choisissons~$n>n_0$ tel que 
$$\mathrm{H}^p({Y}_1,\: \mathrm{R}^q{f}_{1,*}\: \mathcal{K}_Z\otimes \mathcal{L}^n) \: = \: 0$$
pour tout~$p>0$ et~$q\geq 0$ ; c'est possible d'après le critère d'amplitude de Serre. Par la suite spectrale de Leray, on a donc
$$\mathrm{H}^0({Y}_1,\: \mathrm{R}^q{f}_{1,*} \: \mathcal{K}_Z\otimes \mathcal{L}^n) \: = \: \mathrm{H}^q({Z}_1,\: \mathcal{K}_Z\otimes f^* \mathcal{L}^n)\: = \: 0$$
pour tout~$q\geq 0$. Choisissons~$n$ encore plus grand pour que $\mathrm{R}^q{f}_{1,*} \: \mathcal{K}_Z\otimes \mathcal{L}^n$ soit engendré par ses sections globales sur~${Y}_1$ pour tout~$q>0$. On en déduit que $\mathrm{R}^q{f}_{1,*} \: \mathcal{K}_Z\otimes \mathcal{L}^n=0$ donc que $\mathrm{R}^q {f}_{1,*} ( \mathcal{K}_Z) =0$.
\end{demo}

On en déduit par un procédé habituel que pour tout corps~$K$ de caractéristique nulle, tout morphisme $f:Z\rightarrow Y$ birationnel entre schémas propres sur~$\Spec(K)$ et tout diviseur à croisement normaux~$D\subset Z$ d'opposé ample pour~$f$, on a~$\mathrm{R}^q {f}_{*}\:  (\mathcal{K}_{Z}) = 0$ pour tout $q>\nolinebreak 0$. Ce dernier énoncé est un cas particulier du théorème d'annulation relative de Grauert-Riemenschneider~\cite{GR}, qui prédit que pour tout morphisme $f:Z\rightarrow Y$ birationnel entre schémas propres sur~$\Spec(\C)$, on a~$\mathrm{R}^q {f}_{*}\:  (\mathcal{K}_{Z}) = 0$ pour tout~$q>0$. Il n'est pas clair que le théorème de Grauert-Riemenschneider reste vrai en caractéristique positive, même assez grande et même si les objets se relèvent aux seconds vecteurs de Witt~\cite{EV}.

\medskip

\begin{demoth1} Des dévissages habituels permettent de supposer que~$k$ est de caractéristique positive. Il existe alors par hypothèse des objets~$X$, $X^*$, $\bar{X}$, $D$ et~$\pi$ sur~$W_2(k)$.  La proposition~\ref{prop_generale} montre que $\mathrm{R}^q\pi_* \mathcal{K}_{\bar{X}}=0$ sur~$X^*\times \Spec(k)$ pour tout~$q>0$. Mais~$\mathcal{K}_{\bar{X}}=(\pi^*\omega^{\otimes n})(-D)$ par l'isomorphisme de Kodaira-Spencer. On trouve par la formule de projection que~$\omega^n\otimes \mathrm{R}^q\pi_* \mathcal{O}_{\bar{X}}(-D)$ est nul sur~$X^*\times \Spec(k)$ pour tout~$q>0$ et il en est  de même pour~$\mathrm{R}^q\pi_* \mathcal{O}_{\bar{X}}(-D)$.
\end{demoth1}

\providecommand{\bysame}{\leavevmode ---\ }
\providecommand{\og}{}
\providecommand{\fg}{}
\providecommand{\smfandname}{et}
\providecommand{\smfedsname}{\'eds.}
\providecommand{\smfedname}{\'ed.}
\providecommand{\smfmastersthesisname}{M\'emoire}
\providecommand{\smfphdthesisname}{Th\`ese}

\end{document}